# RECOVERING ASYMPTOTICS OF METRICS FROM FIXED ENERGY SCATTERING DATA

MARK S. JOSHI AND ANTÔNIO SÁ BARRETO

ABSTRACT. The problem of recovering the asymptotics of a short range perturbation of the Euclidean metric on $\mathbb{R}^n$ from fixed energy scattering data is studied. It is shown that if two such metrics, $g_1, g_2$, have scattering data at some fixed energy which are equal up to smoothing, then there exists a diffeomorphism $\psi$ 'fixing infinity' such that $\psi^* g_1 - g_2$ is rapidly decreasing. Given the scattering matrix at two energies, it is shown that the asymptotics of a metric and a short range potential can be determined simultaneously. These results also hold for a wide class of scattering manifolds.

## 1. INTRODUCTION

In this paper, we examine the question of the recovery of the asymptotics of a metric from fixed energy scattering data for short range perturbations of Euclidean space. We show that, modulo an inevitable diffeomorphism invariance, the asymptotics are determined. This appears to be the first result on the recovery of asymptotics of a metric, even when given the scattering matrix at all energies. Our approach is to use the techniques of [4] and invert the arising integral transforms which are considerably more complicated than those which appear there. We work in the general context of a manifold equipped with scattering metric and only specialize to $\mathbb{R}^n$ at the final stage.

Recall that a scattering metric on a manifold with boundary, $(X, \partial X)$, is a metric which takes the form
$$g = \frac{dx^2}{x^4} + \frac{h}{x^2}$$
for some boundary defining function $x$ with $h$ a symmetric $2$−cotensor restricting to a positive definite metric on the boundary. We recall the definition of the scattering matrix in this context from [6]. Let $\Delta$ be the Laplacian for $g$. For any non-zero real $\lambda$, given a smooth function $f$ on the boundary, there is a unique smooth function $u$ on the interior of $X$, of the form
$$u = e^{i\frac{\lambda}{x}} x^{\frac{n-1}{2}} f' + e^{-i\frac{\lambda}{x}} x^{\frac{n-1}{2}} f'',$$
with $f', f''$ smooth on $X$ up to the boundary and $f$ equal to the restriction of $f'$ to the boundary, such that
$$(\Delta - \lambda^2) u = 0$$

The scattering matrix is then the map
(1.1) $$S(\lambda) : f \longmapsto f''_{|\partial X}.$$

Melrose and Zworski showed in [8] that $S(\lambda)$ is a classical Fourier integral operator of order $0$ associated to geodesic flow at time $\pi$.

These definitions generalize those of the scattering matrix on $\mathbb{R}^n$, as $\mathbb{R}^n$ can be compactified via stereographic projection. Let
$$\mathbb{S}^n_+ = \{ y \in \mathbb{R}^{n+1} : |y| = 1, \ y_1 > 0 \},$$

then the projection is given by the map

$$SP : \mathbb{R}^n \longrightarrow \mathbb{S}^n_+$$
$$z \longmapsto \left( \frac{1}{(1+|z|^2)^{\frac{1}{2}}}, \frac{z}{(1+|z|^2)^{\frac{1}{2}}} \right).$$

It is not difficult to see that perturbations of the Euclidean metric of the form

$$(1.2) \qquad g_{ij} = \delta_{ij} + \frac{1}{|z|^2} h_{ij}\left(\frac{1}{|z|}, \frac{z}{|z|}\right), \text{ as } |z| \to \infty,$$

where $h_{ij}$ are smooth functions in $[0,1) \times \mathbb{S}^{n-1}$, push-forward under $SP$ to scattering metrics in $\mathbb{S}^n_+$, as in equation (0.9) of [8]. These shall constitute short-range perturbations in a sense we make precise below.

In this setting the theorem of Melrose and Zworski is equivalent to saying that $S(\lambda) = Pa^*$ with $a^*$ pull-back by the antipodal map and $P$ a zeroth order, classical, pseudo-differential operator. These results remain true when a short range potential of the form $x^2 C^\infty(X)$ is added to $\Delta$. It was shown in [4] that for a fixed metric, the asymptotics of such a potential are recoverable from the scattering matrix at one energy by studying the total symbol of this Fourier integral operator. This was extended to potentials of the form $Cx + \mathcal{O}(x^2)$ in [3] for a suitably modified definition of the scattering matrix.

If $\psi$ is a diffeomorphism of $(X, \partial X)$, fixing $\partial X$ and $dx$ on $\partial X$, then the scattering matrix will be left invariant under pull-back of the metric by $\psi$. So we can not in general expect to recover a metric's asymptotics. We establish a normal form for a scattering metric near $\partial X$. We show that there exists a diffeomorphism from $[0, \epsilon) \times \partial X$ to a neighbourhood $U \subset X$ of $\partial X$ so that the pulled-back metric on $[0, \epsilon) \times \partial X$ is of the form

$$\frac{dx^2}{x^4} + \frac{h(x, y, dy)}{x^2} + \mathcal{O}(x^\infty).$$

The recovery of the metric's asymptotics is then equivalent to recovering the Taylor series of $h$ at $x = 0$.

We recall from [4],

**Definition 1.1.** *We define a short range perturbation of a scattering metric $g$ on a manifold $X$ to be a differential operator $F \in \Psi^{2,2}_{sc}(X, \Omega^{\frac{1}{2}}_{sc})$ such that $\Delta_g + F$ is equal to the Laplacian induced by some metric on $X$ plus a smooth potential.*

Note that although this definition is phrased in terms of operators operating on scattering half-densities, any scattering metric will induce a canonical trivialization of that bundle. So adding a short range potential and/or a short range perturbation of the metric yields a short range perturbation. We remark that, in general, the metric associated to $\Delta_g + F$ is necessarily a scattering metric.

The following Theorem was proved in [4],

**Theorem 1.1.** *Let $(X, \partial X, g, x)$ be a smooth manifold with boundary with scattering metric $g$ and boundary defining function $x$. Suppose $F_1, F_2$ are short range perturbations and $F = F_1 - F_2 \in \Psi^{2,k}_{sc}(X, \Omega^{\frac{1}{2}}_{sc})$ where $k \in \mathbb{N}$ is greater than or equal to 2. Let $f(\tau, |\mu|, y, \hat{\mu})$ be the principal symbol of $F$ at the boundary. Let $S_j(\lambda)$ be the scattering matrix associated to $\Delta + F_j - \lambda^2$. Then*

$$S_1(\lambda) - S_2(\lambda) \in I^{-k+1}(G_\pi)$$

*and its principal symbol determines and is determined by*

$$I_{k-1}(f, \gamma) = \int_0^\pi f(|\lambda| \cos s, |\lambda| \sin s, \gamma(s))(\sin s)^{k-1} ds$$

*for every integral curve, $\gamma$, of the geodesic flow on the unit cosphere bundle over the boundary.*

We deduce from this that given two scattering metrics which are equal up to some order at $x = 0$, then we can recover a weighted integral of the next term in the Taylor series along lifted geodesics in the cosphere bundle of the boundary. Such directionally dependent X-ray transforms have previously arisen as the linearization of the problem of recovery of a metric from its hodograph - that is given the length of geodesics do they determine a manifold? Our theorem then follows from judiciously applying some results of Michel, [9].

We shall say that an open set $U$ in $\mathbb{R}^n_y$ is a neighbourhood of infinity if for some $R > 0$, it contains the set $|y| \geq R$. Our main result is,

**Theorem 1.2.** *Suppose $g_1, g_2$ are short range perturbations of the Euclidean metric on $\mathbb{R}^n$, i.e satisfy (1.2), for $n \geq 3$, and their associated scattering matrices are equal up to smoothing at some non-zero energy. Then there exist neighbourhoods of infinity $U, V$ and a diffeomorphism $\psi : U \to V$ which takes the form*
$$\psi(r, \omega) = (r, \omega) + (r^{-1}f(r, \omega), r^{-1}h(r, \omega)),$$
*in polar coordinates, with $f, h$ bounded smooth functions, such that*
$$\psi^* g_1 - g_2 = \mathcal{O}(r^{-\infty}).$$

We remark that $\psi$ will in fact be more constrained than this, as it will be smooth up to the sphere at infinity - this constraint is best understood in the context of scattering on a manifold with boundary, as below. We also deduce, because of the different homogeneities involved in metric and potential scattering, that given the scattering matrix at two energies we can recover the asymptotics of both a metric and potential.

Our theorem holds for certain classes of scattering manifolds other than $\mathbb{R}^n$. It is only the geometry of the boundary that is important and applying our results and those of Michel, [9], we have

**Theorem 1.3.** *Let $(X, \partial X)$ be a manifold with boundary and of dimension $n \geq 3$ with scattering metric $g$, and suppose that each connected component of $\partial X$, with the induced metric, is one of the following,*
- *a sphere of irrational radius,*
- *a sphere of radius one,*
- *real projective space,*

*then if $g_1, g_2$ are short range perturbations of $g$ which have the same scattering matrix at some non-zero energy, up to smoothing, then there exist $U, V$ open neighbourhoods of the boundary and a diffeomorphism, $\psi : U \to V$, fixing the boundary such that*
$$\psi^* g_1 - g_2 = \mathcal{O}(x^\infty)$$
*and $\psi^* x = x + \mathcal{O}(x^3)$ for any boundary defining function $x$.*

The problem of recovering a metric from the scattering matrix has been extensively studied in recent years. Working with the Dirichlet to Neumann map for the wave equation, which can be recovered from the scattering matrix at all energies, see for example section 3.7 of [7], Belishev and Kurylev, [1] and [2], have shown that that, modulo a group of diffeomorphism, smooth compactly supported metric perturbation can be recovered from the knowledge of $S(\lambda)$ for all $\lambda \in \mathbb{R}^n$. In the fixed energy case much less is known. For the related problem of recovering a metric on a compact smooth domain from the Dirichlet to Neumann map for Laplace's equation, Lee and Uhlmann [5] have shown that the Taylor series at the boundary is determined. They use this result to show that, under certain additional assumptions, real analytic metrics are determined, modulo a group of diffeomorphisms, by the Dirichlet to Neumann map for Laplace's equation. In the case of a compactly supported perturbation, the scattering matrix is equal to the identity plus a smoothing operator, so no information can be extracted from its symbol.

The problem of determining a potential from the scattering matrix has a longer history, see [4] for a brief account. As it is mentioned there, there are non-zero potentials in $\mathbb{R}^n$ whose scattering matrix at a fixed energy is the identity. These are called transparent potentials. The question whether there are transparent metrics seems to be an interesting one.

This paper was prepared whilst visiting the Fields Institute and we would like to thank them for their hospitality. We would like to thank Harold Donnely, Richard Melrose, Todd Quinto, Steve Zelditch and Maciej Zworski for helpful comments. We are grateful to Gunther Uhlmann for referring us to Michel's paper [9]. The work of second author was partly supported by the NSF grant DMS-9623175.

## 2. A Normal Form for Scattering Metrics

In this section we establish the existence of boundary normal coordinates for a scattering metric up to infinite order terms. The main result is the following

**Proposition 2.1.** *Let $M$ be a smooth compact manifold with boundary $\partial M$. Let $x \in C^\infty(M)$ be a defining function of $\partial M$ and let*

$$(2.1) \qquad g = \frac{dx^2}{x^4} + \frac{h}{x^2}$$

*be a scattering metric in $M$. Then there exist a neighbourhood $U \subset M$ of $\partial M$ and a diffeomorphism*

$$\Psi : [0, \epsilon) \times \partial M \longrightarrow U, \quad \epsilon > 0,$$

*such that*

1. $\Psi(0, Y) = Y$,
2. $x \circ \Psi(X, Y) = X + X^3 F$, $F \in C^\infty([0, \epsilon) \times \partial M)$,
3. $\Psi^* g = \frac{dX^2}{X^4} + \frac{h'}{X^2} + \mathcal{O}(X^\infty)$, *where $h'(\partial_X, .) = 0$.*

*Proof.* This is the analogue, modulo $\mathcal{O}(X^\infty)$ terms, of boundary normal coordinates on a Riemannian manifold. The difficulty here is, as $g$ is singular at $\partial M$, we do not have an obvious notion of geodesics which are normal to $\partial M$. It is possible that there exists a map $\Psi$ for which 3 holds exactly.

First, for convenience, we identify a neighbourhood $U \subset M$ of $\partial M$ and $[0, \epsilon) \times \partial M$. This can be achieved in the following way. Let $\nabla_x$ be the vector field defined by $g(\nabla_x, \cdot) = dx$ and denote $||\nabla_x||^2 = g(\nabla_x, \nabla_x)$. It follows from (2.1) that the vector field $\Sigma = \frac{1}{||\nabla_x||^2} \nabla_x$ is a smooth vector field in $M$ and moreover $\Sigma x = 1$. For $p \in \partial M$ let $\gamma_\Sigma(p, s)$ denote the point on the integral curve of $\Sigma$ at time $s$ passing through $p$ at time 0. Let $\chi : [0, \epsilon) \times \partial M \longmapsto U$ be the map defined by $\chi(s, p) = \gamma_\Sigma(p, s)$. Since $\Sigma x = 1$, we deduce that $x(\gamma_\Sigma(s, p)) = s$.

Let $p \in \partial M$ and let $y_j$, $1 \leq j \leq n-1$, be smooth functions on $\partial M$ which form a coordinate system in a neighbourhood of $p$. Extending $y_j$ to be constant along the integral curves of $\Sigma$ gives a coordinate system in a neighbourhood of $(0, p)$ in $[0, \epsilon) \times \partial M$.

It follows from (2.1) that, in these coordinates, there exists $r \in \mathbb{N}$, $r \geq 0$, such that $g$ is given by

$$(2.2) \qquad \chi^* g = (1 + x^{2+r} f_0(x, y)) \frac{dx^2}{x^4} + \frac{x^r}{x^2} \sum_{j=1}^{n-1} h_j(x, y) dx dy_j + \frac{1}{x^2} \sum_{i,j=1}^{n-1} h_{ij}(x, y) dy_i dy_j,$$

For simplicity, from now on we use $g$ to denote $\chi^* g$. We want to construct coordinates $(x, y)$ in which $f_0$ and $h_j$, $1 \leq j \leq n-1$, vanish to infinite order at $x = 0$. To do that we in principle do the following. Set

$$(2.3) \qquad x = X + X^{2+r} F(X, Y), \quad y = Y + X^{1+r} G(X, Y),$$

then substitute these expressions into (2.2). The goal is to construct the coefficients of the Taylor series of $F$ and $G$ at $X = 0$, and use Borel's lemma to construct the map (2.3) such that, in terms of $(X, Y)$, $g$ satisfies 3. This gives a rather complicated system of non-linear differential equations satisfied by $F$ and $G$. By analyzing this system we can see that once $F(0, Y)$ is chosen, then the Taylor series of $F$ and $G$ at $X = 0$ are determined. However, to prove this directly would generate a long computation. Instead we shall proceed by induction in $r$.

We set

(2.4) $$x = X + X^{3+r}F(Y), \quad y_j = Y_j + X^{1+r}G_j(Y).$$

We want to choose $F$ and $G_j$ such that, in coordinates $(X, Y)$, $g$ satisfies (2.2) with $r$ replaced by $r + 1$. We remark that this corresponds to setting $F(0, Y) = 0$ in (2.3). Denote $\Psi_0(X, Y) = (x, y)$, $F_0 = f_0 \circ \Psi_0$, $H_j = h_j \circ \Psi_0$, $H_{ij} = h_{ij} \circ \Psi_0$. Then $\Psi_0^* g$ is given by

$$\Psi^* g = A(X, Y)\frac{dX^2}{X^4} + \frac{1}{X^2}\sum_{k=1}^{n-1} A_k(X, Y)dXdY_k + \frac{1}{X^2}\sum_{k,l=1}^{n-1} A_{kl}(X, Y)dY_k dY_l,$$

where,

(2.5)

$$A(X, Y) = 1 + X^2\left(F_0(0, Y) + 2F(Y) + \sum_{j=1}^{n-1}(H_j G_j)(0, Y) + \sum_{ij=1}^{n-1}(H_{ij}G_i G_j)(0, Y)\right) + \mathcal{O}(X^3) \text{ if } r = 0,$$

$$A(X, Y) = 1 + X^{r+2}(F_0(0, Y) + 2(r+1)F(Y)) + \mathcal{O}(X^{r+3}), \text{ if } r > 0,$$

and, using the symmetry of $H_{ij}$, we obtain

(2.6) $$A_k(X, Y) = X^r\left(H_k(0, Y) + 2(r+1)\sum_{j=1}^{n-1} H_{kj}(0, Y)G_j(Y)\right) + \mathcal{O}(X^{r+1}).$$

Since, by assumption, $H_{kj}(0, Y)$ is invertible, we can find $G_j(Y)$, $1 \leq j \leq n-1$, for which $A_k(X, Y) = \mathcal{O}(X^{r+1})$. If $r = 0$, substitute $G_j$ obtained from (2.6) into the first equation in (2.5) and find $F(Y)$ such that $A(X, Y) = \mathcal{O}(X^3)$. If $r > 0$, take $F(Y) = \frac{-1}{2(1+r)}F_0(0, Y)$. Thus, in coordinates $(X, Y)$, $g$ satisfies (2.2) with $r$ replaced by $r + 1$. We observe that, by construction, $F$ and $G_j$ are determined by the metric at $\partial M$ and therefore are smooth functions defined on $\partial M$. This defines a map $\Psi_0 : [0, \epsilon_0) \times \partial M \longmapsto [0, \delta_0) \times \partial M$

Repeating this construction gives a sequence of maps,

$$\phi_k = \Psi_k^{-1} : [0, \delta_k) \times \partial M \longmapsto [0, \epsilon_k) \times \partial M,$$

$k \in \mathbb{N}$ such that

(2.7) $$\phi_k(x, y) = \text{Id} + (x^{r+k+3}F_k, x^{r+k+1}G_k), \quad k \in \mathbb{N}, \ F_k, G_k \in C^\infty([0, \delta_k) \times \partial M),$$

and that $\Phi_j = \phi_j \circ \phi_{j-1} \circ ... \circ \phi_0$, is such that $(\Phi_j^{-1})^* g$ satisfies (2.2) with $r$ replaced by $r + j + 1$.

From (2.7) we deduce that

$$\Phi_{j+l} - \Phi_j = (\mathcal{O}(x^{j+r+3}), \mathcal{O}(x^{j+r+1})), \quad \forall l \geq 0.$$

Thus the coefficients of the Taylor series of the components of $\Phi_{j+l}$, $l \geq 0$, up to order $j + r + 2$ for the first, and $j + r$ for the others, are determined by $\Phi_j$.

Although this construction has been local in $y$, the Taylor series is determined uniquely and so different local coordinates will not change it. There is therefore no problem in passing to the whole of $\partial M$. Hence we can use Borel's lemma to construct a smooth map $\Phi : [0, \epsilon) \times \partial M \longmapsto [0, \epsilon) \times \partial M$ of the form

$$\Phi = \text{Id} + (x^{r+3}F, x^{r+1}G),$$

such that

$$\Phi - \Phi_j = (\mathcal{O}(x^{j+r+3}), \mathcal{O}(x^{j+r+1})),$$

It is then easy to see that $\Psi = \chi \circ \Phi^{-1}$ satisfies properties 1, 2 and 3. This concludes the proof of the Proposition. $\square$

*Remark 1.* The Taylor series of the diffeomorphism is determined once the boundary defining function has been chosen. However, choosing a different boundary defining of the form $X = x + \alpha(y)x^2$ will give rise to a different Taylor series.

## 3. Proof of the Theorem

In the previous section we have established that given two scattering metrics, $g_l$, $l = 1, 2$, there are diffeomorphisms $\Psi_l$ from neighbourhoods of the boundary to $[0, \epsilon) \times \partial X$, so that $\Psi_l^* g_l$ is of the form

$$\frac{dx^2}{x^4} + \frac{1}{x^2} h_l + \mathcal{O}(x^\infty), \text{ with } h_l(\partial_x, \cdot) = 0. \tag{3.1}$$

As we are working modulo diffeomorphism invariance, we can therefore henceforth assume that we are working in a product decomposition in which both metrics have this form.

Theorem 1.3 is an immediate consequence of the following

**Proposition 3.1.** *Let $\Psi_l^* g_l$ be of the form (3.1). Suppose that each component of $\partial X$ is either a sphere of radius one, the real projective space or a sphere of irrational radius. Let $h_0, h_1$ and $h_{j,l}$, $1 \leq j \leq k$, $l = 1, 2$, be symmetric 2-cotensors defined on $\partial X$. Suppose that*

$$h_l = h_0 + x h_1 + \tilde{h}_l + x^k h_{k,l} + x^{k+1} R_l, \tag{3.2}$$

*where $\tilde{h}_l = \sum_{j=2}^{k-1} x^j h_{j,l}$. If for a fixed non-zero, energy $\lambda$,*

$$S_1(\lambda) - S_2(\lambda) \in I^{-k}(G_\pi)$$

*then $h_{j,1} = h_{j,2}$ for $j = 1, ..., k$.*

As there is no interaction between the different components of the boundary at the level of singularities of the scattering matrix, we can deal with each of the components individually. Our goal is, of course, to apply Theorem 1.1 to deduce that a weighted integral transform of $h_{k,1} - h_{k,2}$ is zero and then show that this implies that $h_{k,1} - h_{k,2}$ is zero, too. We remark that, since we assume that the perturbations are short range, and $k \geq 2$ in Theorem 1.1, we need to impose that the first two terms in the Taylor's expansion are independent of $l$. To achieve our goal we first need to compute the principal symbol of the $k$-th order perturbation at the boundary, $k \geq 2$. This is done in the following

**Lemma 3.1.** *Let $X$ be a $C^\infty$ compact manifold with boundary, $\partial X$, of dimension $n$, let $g$ be a scattering metric on $X$. Let $x$ be a defining function of $\partial X$ and let $(x, y)$, $y = (y_1, ..., y_{n-1})$ be local coordinates valid near $q_0 \in \partial M$ in which $g$ is of the form*

$$g = \frac{dx^2}{x^4} + \frac{1}{x^2} h(x, y, dy) + \mathcal{O}(x^\infty)$$

$$h(x, y, dy) = h_0(y, dy) + \tilde{h}(x, y, dy) + x^k h_k(y, dy) + x^{k+1} R(x, y, dy), \tag{3.3}$$

$$\tilde{h}(x, y, dy) = \sum_{m=1}^{k-1} x^m h_m(y, dy)$$

*Let $H = (h_{ij})$, $H_0 = \left((h_0 + \tilde{h})_{ij}\right)$ and $H_k = ((h_k)_{ij})$, be the matrices of coefficients of $h$, $h_0 + \tilde{h}$ and $h_k$ respectively. Let $\delta_0 = \det H_0$, $H_0^{-1} = \left((h_0 + \tilde{h})^{ij}\right)$ and $B_k = \left(B_k^{ij}\right) = H_0^{-1} H_k H_0^{-1}$. Let $\Delta_g$ denote the Laplacian with respect to $g$. Then*

(3.4)
$$\Delta_g = \frac{x^{n+1}}{(\det H)^{\frac{1}{2}}}\partial_x \frac{(\det H)^{\frac{1}{2}}}{x^{n+1}}x^4\partial_x + \frac{x^2}{\delta_0^{\frac{1}{2}}}\sum_{i,j=1}^{n-1}\partial_{y_i}\delta_0^{\frac{1}{2}}h_0^{ij}\partial_{y_j} - \frac{x^{k+2}}{\delta_0^{\frac{1}{2}}}\sum_{i,j=1}^{n-1}\partial_{y_i}\delta_0^{\frac{1}{2}}B_k^{ij}\partial_{y_j}$$
$$+\frac{x^{k+2}}{2\delta_0^{\frac{1}{2}}}\sum_{i,j=1}^{n-1}D_{y_i}\delta_0^{\frac{1}{2}}\operatorname{Tr}(H_0^{-1}H_k)h_0^{ij}\partial_{y_j} - \frac{x^{k+2}}{2\delta_0^{\frac{1}{2}}}\operatorname{Tr}(H_0^{-1}H_k)\sum_{i,j=1}^{n-1}\partial_{y_i}\delta_0^{\frac{1}{2}}h_0^{ij}\partial_{y_j} + x^{k+3}Q(x,y,\partial_y) + \mathcal{O}(x^\infty),$$

where $Q$ is a smooth second order differential operator and $\mathcal{O}(x^\infty)$ denotes a second order operator whose coefficients vanish to infinite order at $\partial M$.

*Proof.* Let $(g_{ij}) = G$, $(g^{ij}) = G^{-1}$, $\delta = \det G$, and denote $x_0 = x$, $x_j = y_j$, $1 \leq j \leq n-1$. Then, by definition, $\Delta_g = \frac{1}{\delta^{\frac{1}{2}}}\sum_{i,j=0}^{n}\partial_{x_i}g^{ij}\delta^{\frac{1}{2}}\partial_{x_j}$. We know that $g_{00} = \frac{1}{x^4}$, $g_{0j} = 0$, $g_{ij} = \frac{1}{x^2}h_{ij}$, $1 \leq i,j \leq n-1$. Hence the inverse matrix is given by $g^{00} = x^4$, $g^{0j} = 0$, $g^{ij} = x^2 h^{ij}$, $1 \leq i,j \leq n-1$. From (3.3) we obtain,

(3.5) $$H = H_0(I + x^k H_0^{-1}H_k + x^{k+1}R_1),$$

where here, and in what follows, $R_1$ denotes a smooth matrix. Therefore

(3.6) $$H^{-1} = (I + x^k H_0^{-1}H_k + x^{k+1}R_1)^{-1}H_0^{-1} = H_0^{-1} - x^k H_0^{-1}H_k H_0^{-1} + x^{k+1}R_1.$$

We know, from the definition of $g$, that

(3.7) $$\det H = (\det H_0)\det(I + x^k H_0^{-1}H_k + x^{k+1}R_1) = (\det H_0)(1 + x^k \operatorname{Tr}(H_0^{-1}H_k) + x^{k+1}F_1),$$

where $F_1$ is a smooth function. From the definition of $g$ and $\Delta_g$ we obtain

(3.8) $$\Delta_g = \frac{x^{n+1}}{(\det H)^{\frac{1}{2}}}\partial_x \frac{(\det H)^{\frac{1}{2}}}{x^{n+1}}x^4\partial_x + \frac{x^2}{(\det H)^{\frac{1}{2}}}\sum_{i,j=1}^{n-1}\partial_{y_i}\frac{(\det H)^{\frac{1}{2}}}{x^{n+1}}h^{ij}\partial_{y_j}.$$

From (3.6) we know that, for $B_k = H_0^{-1}H_k H_0^{-1}$,
$$h^{ij} = h_0^{ij} - x^k B_k^{ij} + x^{k+1}R^{ij}.$$

We deduce from (3.7) and (3.8) that

(3.9)
$$\Delta_g - x^{k+3}Q(x,y,\partial_y) - \mathcal{O}(x^\infty) = \frac{x^{n+1}}{(\det H)^{\frac{1}{2}}}\partial_x\frac{(\det H)^{\frac{1}{2}}}{x^{n+1}}x^4\partial_x + \frac{x^2}{(\det H_0)^{\frac{1}{2}}}\partial_{y_i}(\det H_0)^{\frac{1}{2}}h_0^{ij}\partial_{y_j}$$
$$+\frac{x^{k+2}}{(\det H_0)^{\frac{1}{2}}}\partial_{y_i}(\det H_0)^{\frac{1}{2}}\left(-B_k^{ij} + \frac{1}{2}\operatorname{Tr}(H_0^{-1}H_k)h_0^{ij}\right)\partial_{y_j} - \frac{x^{k+2}}{2(\det H_0)^{\frac{1}{2}}}\operatorname{Tr}(H_0^{-1}H_k)\partial_{y_i}(\det H_0)^{\frac{1}{2}}h_0^{ij}\partial_{y_j}.$$

Now (3.4) follows directly from (3.9). □

The following is then an immediate consequence of Proposition 3.1,

**Corollary 3.1.** *Let $g_l$, $l = 1,2$ be given by (3.1)-(3.2) and equal up to $k^{th}$ order with $k \geq 2$. Then*
$$F = \Delta_{g_1} - \Delta_{g_2} \in \Psi_{sc}^{2,k}([0,\epsilon) \times \partial X, \Omega_{sc}^{\frac{1}{2}}).$$
*Let $f_k$ be the principal symbol of $F$ at $\partial X$. Then*

(3.10) $$f_k(y,\mu) = T_k(y,\mu),$$

where, in the notation of Lemma 3.1, , $T_k$ is the 2-cotensor on $T^*\partial X$, defined in local coordinates by

$$T_k(y,\mu) = \sum_{i,j=1}^{n-1} B_k^{ij}(0,y)\mu_i\mu_j.$$

We remark that, as $H_0(0,y)$ is the matrix of coefficients of $h_0(y)$, which we also denote by $h_0$,

$$B_k(0,y) = h_0^{-1} H_k h_0^{-1}.$$

The reason for introducing coordinates in which $g_l$, $l = 1, 2$, is given by (3.1) is that $f_k$ does not depend on $\tau$, the coefficient of $\frac{dx}{x^2}$. So far we have obtained the integrand of the transformation in Theorem 1.1. Next we show that this transformation, for this type of integrand, is injective, when $\partial X$ is one of the manifolds in Theorem 1.3.

**Lemma 3.2.** *Let $F$ be a smooth symmetric 2-cotensor on $\mathbb{S}^{n-1}$. We denote $F(x,V)$ its value at the point $x \in \mathbb{S}^{n-1}$ and the $V \in T_x\mathbb{S}^{n-1}$. If there exists $k \in \mathbb{N}$ such that for every geodesic $\gamma$ of length $\pi$, with respect to the Euclidean metric,*

$$(3.11) \qquad I_k(F,\gamma) = \int_0^\pi F(\gamma(t),\gamma'(t),\gamma'(t))(\sin t)^k dt = 0$$

*Then*
*i) $F$ has the parity of $k$, i.e*

$$F(\gamma(t),\gamma'(t),\gamma'(t)) = (-1)^k F(\gamma(t+\pi),\gamma'(t+\pi),\gamma'(t+\pi))$$

*ii)*

$$(3.12) \qquad \int_0^{2\pi} p(\gamma(t))F(\gamma(t),\gamma'(t),\gamma'(t))dt = 0,$$

*for every polynomial of degree $k$.*

*Proof.* To prove i) we fix a geodesic $\gamma$ and set, for $\alpha \in \mathbb{R}$

$$I_k(\alpha) = \int_0^\pi F(\gamma(t+\alpha),\gamma'(t+\alpha),\gamma'(t+\alpha))(\sin t)^k dt.$$

Since (3.11) holds for every geodesic, we have $I_k(\alpha) = 0$ for all $\alpha$. On the other hand, differentiation with respect to $\alpha$ and integration by parts gives, as in [4],

$$(3.13) \qquad I_k''(\alpha) + k I_k(\alpha) = k(k-1) I_{k-2}.$$

By induction we obtain $I_0(\alpha) = 0$ for $k$ even and $I_1(\alpha) = 0$ for $k$ odd. However, differentiating $I_0$ and $I_1$, we obtain

$$I_0'(\alpha) = F(\gamma(\alpha+\pi),\gamma'(\alpha+\pi),\gamma'(\alpha+\pi)) - F(\gamma(\alpha),\gamma'(\alpha),\gamma'(\alpha)) = 0$$

$$I_1(\alpha) + I_1''(\alpha) = F(\gamma(\alpha+\pi),\gamma'(\alpha+\pi),\gamma'(\alpha+\pi)) + F(\gamma(\alpha),\gamma'(\alpha),\gamma'(\alpha)) = 0.$$

This ends the proof of i). To prove ii), let us assume, for convenience, that $k$ is even. The proof for $k$ odd is the same. We observe that any closed geodesic $\gamma$ in $\mathbb{S}^n$ must intersect the equator $\gamma_e$. Let $x = (x_1, ..., x_{n+1})$ be coordinates in $\mathbb{R}^n$ such that $\gamma_e = \mathbb{S}^n \cap \{x_1 = 0\}$. Now we have that, at time $t$, the $x_1$ coordinate of $\gamma(t)$ is a fixed multiple of $\sin(t)$. Therefore we deduce from (3.11) that (3.12) holds for $p(x) = x_1^k$. By the rotational invariance of $\mathbb{S}^n$, (3.12) must hold for $p(x) = \sum_{j=1}^{n+1} a_j x_j$, with $\sum_{j=1}^{n+1} a_j^2 = 1$. Hence, (3.12) must hold for any homogeneous polynomial of degree $k$. Using (3.13) we deduce the same result for homogeneous polynomials of degree $k-2, k-4, ..., 0$. Thus (3.12) holds for sums of homogeneous polynomials of even degree. On the other hand, from i) we deduce that (3.12) holds for sums of homogeneous polynomials with odd degree. □

Next we state a result of R. Michel, Corollary A' of [10].

**Theorem 3.1.** *Let $F$ be a symmetric 2-cotensor on $\mathbb{S}^{n-1}$ such that, for every closed geodesic $\gamma$ with respect to the Euclidean metric $g_e$,*

$$\int_0^{2\pi} F(\gamma(t), \gamma'(t), \gamma'(t)) dt = 0.$$

*Then there exists a unique decomposition $F = F_1 + F_2$, where $F_1$ is odd and $F_2 = L_X(g_e)$, is the Lie derivative of $g_e$ with respect to a smooth vector field $X$ in $\mathbb{S}^{n-1}$.*

We also recall an equation used in [9] which is satisfied by Lie derivatives of $g_e$. Let

$$\nabla_s : T_1^0(\mathbb{S}^n) \longrightarrow T_2^0(\mathbb{S}^n)$$

be the symmetrized covariant derivative, and let $\nabla_s^*$ be its formal adjoint. Let $X$ be a smooth vector field on $\mathbb{S}^{n-1}$ and $K = L_X(g_e)$, then

(3.14) $$\nabla_s^* \nabla_s^* K + \Delta \operatorname{Tr} K = (n+1)^2 \operatorname{Tr} K.$$

Where $\operatorname{Tr} K$ is the trace of $K$ with respect to $g_e$.

We use these results to prove

**Proposition 3.2.** *Let $F$ be an even symmetric 2-cotensor on $\mathbb{S}^{n-1}$ such that, for every closed geodesic $\gamma$ with respect to the Euclidean metric $g_e$, and for every homogeneous polynomial, $p$, of degree 2,*

$$\int_0^{2\pi} p(\gamma(t)) F(\gamma(t), \gamma'(t), \gamma'(t)) dt = 0,$$

*then $F = 0$.*

*Proof.* From i) of Lemma 3.2 and Theorem 3.1, $p(x)F = K = L_X(g_e)$. Hence it must satisfy (3.14). Recall, see for example [9], that in local coordinates where $K = \sum_{i,j=1}^{n-1} K_{ij} dx_i dx_j$, its trace is given by $\operatorname{Tr}(K) = \sum_{i,j=1}^{n-1} K_{ij} g_e^{ij}$, and

$$(\nabla_s^* K)_m = \sum_{k,l=1}^{n-1} \partial_{x_l} K_{mk} g_e^{kl} \text{ and } \nabla_s^* \nabla_s^* K = \sum_{l,m=1}^{n-1} \partial_{x_l} (\nabla_s^* K)_m g_e^{lm}.$$

Let $n_p$ denote the north pole and let $\mathbb{S}^{n-1}$ be parametrized near $n_p$ by

$$(x_1, ..., x_{n-1}) \mapsto (x_1, x_2, ..., x_{n-1}, (1 - |x|^2)^{\frac{1}{2}})$$

which is valid in the region $x_n > 0$.

In these coordinates the Euclidean metric $g_e$ is given by $g_e = \sum_{j=1}^{n-1} dx_j^2 + G_0$ where $G_0$ vanishes to second order at $n_p$. Hence we deduce that

$$\Delta_{g_e} = \sum_{j=1}^{n-1} \partial_{x_j}^2 + Q(x, \partial_x), \quad \nabla_s^* \nabla_s^* K = \sum_{i,j=1}^{n-1} \partial_{x_i} \partial_{x_j} K_{ij} + Z(x),$$

where $Q$ and $Z$ vanish at $n_p$. Now applying (3.14) to $p(x)F$, for $p(x) = x_j^2$ and $p(x) = x_j x_k$, respectively, we obtain at $n_p$

(3.15) $$2F_{jj} + 2\sum_{i=1}^n F_{ii} = (n+1)^2 \sum_{i=1}^n F_{ii}, \quad 2F_{ij} = (n+1)^2 \sum_{i=1}^n F_{ii}, \quad i \neq j.$$

Adding the first equation in $j$ gives that, at $n_p$, $2\sum_{i=1}^n F_{ii} = n((n+1)^2 - 2)\sum_{i=1}^n F_{ii}$. Hence $\sum_{i=1}^n F_{ii} = 0$, at $n_p$, and from (3.15) we obtain $F_{ij}(n_p) = 0$, $1 \leq i, j \leq n-1$. Since $n_p$ is an arbitrary point on $\mathbb{S}^{n-1}$ we deduce that $F = 0$. This concludes the proof of the Proposition. □

Now we are ready to prove Proposition 3.1.

*Proof.* We proceed by induction. Suppose that for some non-zero $\lambda$, $S_1(\lambda) - S_2(\lambda) \in I^{-2}(G_\pi)$. Then we know from Theorem 1.1 and Corollary 3.1 that

$$\int_0^\pi (T_2^1 - T_2^2)(\gamma(t), \gamma'(t), \gamma'(t))(\sin t)^3 dt = 0,$$

where $T_2^l$ denotes the tensor $T_2$ for the metric $g_l$, $l = 1, 2$. Let us denote $W_2 = T_2^1 - T_2^2$ This guarantees that $W_2(\gamma(t), \gamma'(t), \gamma'(t))$ must vanish for some $t_0 \in [0, \pi]$. Thus, by Lemma 3.2,

$$W_2(\gamma(t_0), \gamma'(t_0), \gamma'(t_0)) = W_2(\gamma(t_0 + k\pi), \gamma'(t_0 + k\pi), \gamma'(t_0 + k\pi)) = 0, \quad k \in \mathbb{Z}.$$

If $\partial X$ is a sphere of irrational radius, the set $\{\gamma(t_0 + k\pi), \ k \in \mathbb{Z}\}$ is dense on $\{\gamma(t), \ t \in [0,1]\}$. Therefore $W_2(\gamma(t), \gamma'(t), \gamma'(t)) = 0$ for all $t \in [0, \pi]$. Since $\gamma$ is arbitrary, this gives that $W_2 = 0$.

When $\partial X = \mathbb{S}^{n-1}$, we use, as in the proof of Lemma 3.2, that along $\gamma(t)$, $x_1$ is a fixed multiple of $\sin(t)$. Hence we deduce that the tensor $\widetilde{W}_2 = x_1(T_2^1 - T_2^2)$ satisfies

$$\int_0^\pi \widetilde{W}_2(\gamma(t), \gamma'(t), \gamma'(t))(\sin t)^2 dt = 0,$$

Thus, from Lemma 3.2, $\widetilde{W}_2$ is even and it must satisfy (3.12) for every polynomial of degree 2. Hence, from Proposition 3.2, $\widetilde{W}_2 = 0$. Hence $W_2 = 0$ and $T_2^1 = T_2^2$. Since $h_0$ is the same for both metrics, we deduce that $h_{2,1} = h_{2,2}$. This shows that the metrics agree to second order at $\partial X$.

Suppose $h_{j,1} = h_{j,2}$ for $1 \leq j \leq k-1$ and that for some $\lambda \neq 0$,

$$S_1(\lambda) - S_2(\lambda) \in I^{-k}(G_\pi).$$

Proceeding as above we deduce that $T_k^1 = T_k^2$. Again using that $h_0$, is the same for both metrics, we deduce that $h_{k,1} = h_{k,2}$. This concludes the proof of the Proposition.

As the space of even functions on the sphere is the same as the space of functions on the projective space, the result when $\partial X$ is projective space is now clear. □

To get the result that the scattering matrix at two energies determines the asymptotics of a short range potential and a short range metric perturbation, we simply observe that at the $k^{th}$ level a metric perturbation will yield a principal symbol in the difference of the scattering matrices which is homogeneous in $\lambda$ of order $k+1$ and a potential perturbation a principal symbol which is homogeneous of order $k-1$. So the two principal symbols are thus determined by the knowledge of the principal symbol for two values of $\lambda$, we can therefore apply the results of [4] and this paper to each individually and the result follows.

Department of Pure Mathematics and Mathematical Statistics,
University of Cambridge
16 Mill Lane,
Cambridge CB2 1SB, England, U.K.
  *E-mail address*: `joshi@@dpmms.cam.ac.uk`

Department of Mathematics
Purdue University,
West Lafayette IN 47907, Indiana, U.S.A.
  *E-mail address*: `sabarre@@math.purdue.edu`